\newfont{\Bbb}{msbm10 scaled\magstephalf}
\documentclass[reqno]{amsart}
\usepackage{amssymb}
\usepackage{amsfonts}
\usepackage{amsmath, amsthm, amscd, amsfonts}

\usepackage{amsmath}
\usepackage{amssymb}
\usepackage{graphicx}%
\providecommand{\U}[1]{\protect\rule{.1in}{.1in}}
\newtheorem{theorem}{Theorem}[section]

\newtheorem{corollary}[theorem]{Corollary}

\newtheorem{definition}[theorem]{Definition}
\newtheorem{example}[theorem]{Example}

\newtheorem{lemma}[theorem]{Lemma}

\newtheorem{proposition}[theorem]{Proposition}
\newtheorem{remark}[theorem]{Remark}

\begin{document}


\title[On Embedding problem of linear fractional maps]{On Embedding problem of linear fractional maps on the unit ball of
$\mathbb{C}^{N}$}

\author[R.Y. Chen and Z.H. Zhou] {Ren-Yu Chen and Ze-Hua Zhou$^*$}
\address{\newline  Ren-Yu Chen, Department of Mathematics,
Tianjin University, Tianjin 300072, P.R. China.}
\email{chenry@tju.edu.cn}
\address{\newline  Ze-Hua Zhou, Department of Mathematics, Tianjin University, Tianjin
300072, P.R. China.
\newline
Center for Applied Mathematics, Tianjin University, Tianjin 300072,
P.R. China.} \email{zehuazhoumath@aliyun.com;
zhzhou@tju.edu.cn}






\date{}
\thanks{\noindent $^*$Corresponding author.\\
This work was supported in part by the National Natural Science
Foundation of China (Grant Nos. 11401426; 11371276; 11301373).}

\begin{abstract}
This paper focuses on the embedding problem of linear fractional maps which
explains when a linear fractional self-map of $B_{N}$ can be a member of a
semigroup of holomorphic self-maps on the unit ball $B_{N}$ of the complex $N$-dimensional Euclidean
space $\mathbb{C}^{N}$.

\end{abstract}
\maketitle


\section{Introduction}

Throughout this article ${B_{N}}$ stands for the open unit ball of the complex $N$-dimensional Euclidean
space $\mathbb{C}^{N}$. A holomorphic mapping $\varphi$ is said to be linear
fractional if%
\[
\varphi\left(  z\right)  =\frac{Az+B}{\left\langle z,C\right\rangle +D},
\]
where $A\in\mathbb{C}^{N\times N}$, $B$ $\in\mathbb{C}^{N}$, $C\in
\mathbb{C}^{N}$ and $D\in\mathbb{C}$. Let $\operatorname{LFM}\left(
B_{N}\right)  $ be the family of linear fractional self-maps on $B_{N}$.
Properties of linear fractional maps have been deeply studied. We refer the
readers to \cite{B1,B2,hyperbolic,BB,JFA1,kms2011,CM1,CM2,JFA2,JFA3,JO} etc.

For each $A\in\mathbb{C}^{N\times N}$, as usual, $\sigma\left(  A\right)  $ is
the spectrum of $A$, $\rho\left(  A\right)  $ is the spectra radius of $A$,
$A^{H}$ is the conjugate transpose of $A$, and $\left\Vert A\right\Vert $ is
the spectra norm of $A$. Let $\exp\left(  A\right)  =\sum_{n=0}^{\infty}%
\frac{A^{n}}{n!}$ be the exponential of $A$. We say that $A$ is dissipative if
for any $w\in\mathbb{C}^{N}$,%
\[
\operatorname{Re}w^{H}Aw\leq0.
\]

The following theorem which was first proved in \cite{Her} is key to the
classification of holomorphic self-maps of $B_{N}$.

\begin{theorem}
[Denjoy-Wolff Theorem on ${B_{N}}$]Let $\varphi$ be a holomorphic self-map of
$B_{N}$. If $\varphi$ has no interior fixed points, then there is a unique
point $w\in\partial B_{N}$ such that the iteration $\left\{  \varphi
_{n}=\underset{n\text{ times}}{\underbrace{\varphi\circ\cdots\circ\varphi}%
}\right\}  $ converges to $w$ uniformly on compact subsets of $B_{N}$.
\end{theorem}

The point $w\in\partial B_{N}$ in the above theorem is called the Denjoy-Wolff
point of $\varphi$. According to Theorem 1.3 in \cite{Ma}, there is a real
number $\delta\in(0,1]$ such that
\[
\underset{z\rightarrow w}{\lim\inf}\dfrac{1-|\varphi(z)|^{2}}{1-|z|^{2}%
}=\delta.
\]
In this case, $\delta$ is said to be the boundary dilation coefficient of
$\varphi$. With the help of the above definitions, the holomorphic self-maps
of $B_{N}$ can be classified into three groups:

\begin{definition}
Let $\varphi$ be a holomorphic self-map of $B_{N}$.

\begin{enumerate}
\item If $\varphi$ has at least one interior fixed point then it is called elliptic;

\item If $\varphi$ has no interior fixed point and its boundary dilation
coefficient $\delta\in\left(  0,1\right)  $, then $\varphi$ is called hyperbolic;

\item If $\varphi$ has no interior fixed point and boundary dilation
coefficient $\delta=1$, then $\varphi$ is called parabolic.
\end{enumerate}
\end{definition}

A continuous semigroup $\left\{  \varphi_{t}\right\}  $ of holomorphic
mappings on a domain $D$ of $\mathbb{C}^{N}$ is a continuous homomorphism from
the additive semigroup of non-negative real numbers into the composition
semigroup of all holomorphic self-maps of $D$ endowed with the compact-open
topology. To each continuous semigroup $\left\{  \varphi_{t}\right\}  $ there
corresponds exactly a holomorphic vector field $F:D\rightarrow\mathbb{C}^{N}$
such that%
\[
\frac{\partial\varphi_{t}}{\partial t}=F\circ\varphi_{t}.
\]
The vector field $F$ is called the infinitesimal generator of the semigroup.
For the theory of continuous semigroups we refer to books by Engel and Nagel
\cite{EN1,EN2} and Shoikhet \cite{S}.

An element of a semigroup $\left\{  \varphi_{t}\right\}  $ is said to be an
iterate of $\left\{  \varphi_{t}\right\}  $. There are many special properties
about semigroups, for instance:

\begin{itemize}
\item every iterate of $\left\{  \varphi_{t}\right\}  $ is an injection;

\item if one of the iterates is an automorphism, then all of the iterates are automorphisms;

\item for any $z\in D$, the map $t\mapsto\varphi_{t}\left(  z\right)  $ is
real analytic.
\end{itemize}

In \cite{BCM1}, Bracci, Contreras and D\'{l}az-Madrigal showed that all
elliptic semigroups could be linearized. They provided a basic
\textquotedblleft model\textquotedblright\ for a linear fractional map with no
fixed points in $B_{N}$ (see Theorem \ref{semi4.1}) and then provided a
complete classification up to conjugation of continuous semigroups of linear
fractional self-maps on $B_{N}.$ The following theorem is one of the
conclusions made in \cite{BCM1}.

\begin{theorem}
Let $\left\{  \varphi_{t}\right\}  $ be a semigroup on $B_{N}$. If there is a
$t_{0}\in\left(  0,+\infty\right)  $ such that $\varphi_{t_{0}}$ is an
elliptic (a hyperbolic\ or a parabolic) self-map, then for any $t\in\left(
0,+\infty\right)  $, $\varphi_{t}$ is elliptic (hyperbolic or parabolic).
Furthermore, if $\varphi_{t_{0}}$ is non-elliptic, then all the iterates of
$\left\{  \varphi_{t}\right\}  $ share the same Denjoy-Wolff point.
\end{theorem}

Some other related results about semigroups can be found, for example, in
\cite{AMPA1992}, \cite{JFA6}, \cite{BerAndPor}, \cite{JEMS2010},
\cite{JFA5},\cite{PJM2005}, \cite{JFA4}, and \cite{proc2011}.

A very important problem in the theory of semigroups is that of embedding a
given holomorphic self-map into a semigroup of holomorphic self-maps. We refer
\cite{S} and \cite{FS} for some related results. In \cite{BCM2}, the authors
gave a complete description of infinitesimal generators associated with
semigroup of linear fractional maps on the unit ball of $\mathbb{C}^{N}$. For
the case $N=1$ they showed that a generic semigroup of holomorphic self-maps
of the unit disc is a semigroup of linear fractional maps if and only if it
contains a linear fractional map for some positive time. They also completely
described the associated Koenigs function and solved the embedding problem
from a dynamical point of view. The following theorem is the main result in
\cite{BCM2}.

\begin{theorem}
[{\cite[Theorem 3.3]{BCM2}}]Let $\varphi$ be an arbitrary linear fractional
map of the unit disc $\mathbb{D\subset C}$.

\begin{enumerate}
\item If $\varphi$ is trivial, neutral-elliptic, hyperbolic or parabolic, then
$\varphi$ can be always embedded into a semigroup in $\mathbb{D.}$

\item If $\varphi$ is attractive elliptic with Denjoy-Wolff point $\tau
\in\mathbb{D}$ and repulsive fixed point $\beta\in\mathbb{C}_{\infty
}\backslash\mathbb{D}$, let $\lambda$ be the length of canonical spiral
associated to $\varphi^{\prime}\left(  \tau\right)  \in\mathbb{D}_{{}%
}\backslash\left\{  0\right\}  $. Then $\varphi$ can be embedded into a
semigroup in $\mathbb{D}$ if and only if%
\[
\left\vert \bar{\tau}-\frac{1}{\beta}\right\vert \leq\left\vert \varphi
^{\prime}\left(  \tau\right)  \right\vert \left\vert 1-\frac{\tau}{\beta
}\right\vert .
\]

\end{enumerate}
\end{theorem}

Inspired by the above result, we consider the embedding problem of linear
fractional self-maps on $B_{N}$ in this article.

Let $D\subset\mathbb{C}^{N}$ be a domain and $\psi:D\rightarrow D$ be
holomorphic, $\left(  \psi,D\right)  $ is said to be an iteration couple. Two
couples $\left(  \psi,D\right)  $ and $\left(  \psi^{\prime},D^{\prime
}\right)  $ are said to be conjugated if there exists a biholomorphic map
$\sigma:D\rightarrow D^{\prime}$ such that $\psi=\sigma^{-1}\circ\psi^{\prime
}\circ\sigma$. The map $\sigma$ is called an intertwining map. Sometimes we
will say that $\psi$ and $\psi^{\prime}$ are conjugated instead of $\left(
\psi,D\right)  $ and $\left(  \psi^{\prime},D^{\prime}\right)  $ are
conjugated. Apparently, if $\left(  \psi,D\right)  $ and $\left(  \psi
^{\prime},D^{\prime}\right)  $ is conjugated, then $\psi$ can be embedded into
a semigroup on $D$ if and only if $\psi^{\prime}$ can be embedded into some
other semigroup on $D^{\prime}$.

If $\varphi\in\operatorname{LFM}\left(  B_{N}\right)  $ is elliptic, and
$z_{0}$ is a fixed point of $\varphi$ with $\rho\left(  \varphi^{\prime
}\left(  z_{0}\right)  \right)  =1,$ according to Theorem \ref{ell_normal},
$\left(  \varphi,B_{N}\right)  $ is conjugated to $\left(  \psi,B_{N}\right)
$ where $\psi$ is defined by $\left(  \ref{ell_1}\right)  $. According to
Example \ref{counterexample}, such a $\varphi$ could be embedded into a
semigroup of non-linear fractional maps. The following theorem discusses about
whether such a map can be embedded into a linear fractional semigroup or not.

\begin{theorem}
\label{main1}Let $\varphi\in\operatorname{LFM}\left(  B_{N}\right)  $ be
elliptic. If $\left(  \varphi,B_{N}\right)  $ is conjugated to $\left(
\psi,B_{N}\right)  $, where%
\begin{equation}
\psi\left(  z^{\prime},z^{\prime\prime}\right)  =\left(  \Lambda z^{\prime
},A_{1}z^{\prime\prime}\right)  , \label{ell_1}%
\end{equation}
with $\Lambda$ a diagonal unitary matrix, $\rho\left(  A_{1}\right)  <1$,
$\left\Vert A_{1}\right\Vert \leq1$. Then there is a semigroup $\left\{
\varphi_{t}\right\}  $ of linear fractional maps on $B_{N}$ such that
$\varphi_{1}=\varphi$ if and only if there is a dissipative matrix $M$ with
$\sigma\left(  M\right)  \subset\left\{  \lambda:\operatorname{Re}\lambda
\leq0\right\}  $ such that%
\[
\exp\left(  M\right)  =A_{1}\text{.}%
\]

\end{theorem}

This theorem seems very useful when $N=2$. In fact, it is very easy to prove
the following corollary.

\begin{corollary}
Let $\varphi\in\operatorname{LFM}\left(  B_{2}\right)  $ be elliptic and
$z_{0}\in B_{2}$ be a fixed point of $\varphi$. If $\rho\left(  \varphi
^{\prime}\left(  z_{0}\right)  \right)  =1$, then there is a semigroup
$\left\{  \varphi_{t}\right\}  $ of linear fractional maps of $B_{2}$ such
that $\varphi_{1}=\varphi$.
\end{corollary}

If $\varphi$ $\in\operatorname{LFM}\left(  B_{N}\right)  $ is elliptic with
$z_{0}\in B_{N}$ a fixed point, and $\rho\left(  \varphi^{\prime}\left(
z_{0}\right)  \right)  <1$. Then according to Theorem \ref{ell_normal},
$\left(  \varphi,B_{N}\right)  $ is conjugated to $\left(  \psi,B_{N}\right)
$ where $\psi$ is defined by $\left(  \ref{ell_2}\right)  $. In this case, we
will show in the proof that if $\varphi$ can be embedded into a semigroup,
then this group must be a semigroup of linear fractional maps. And we have the
following theorem.

\begin{theorem}
\label{main2}If $\left(  \varphi,B_{N}\right)  $ is conjugated to $\left(
\psi,B_{N}\right)  $, where%
\begin{equation}
\psi\left(  z\right)  =\frac{Az}{\delta\left\langle z,\left(  A^{H}-E\right)
e_{1}\right\rangle +1}, \label{ell_2}%
\end{equation}
with $A\in\mathbb{C}^{N\times N}$ and $\rho\left(  A\right)  <1$, $\left\Vert
A\right\Vert <1$, $\delta\in\left[  0,1\right]  $. Then there is a semigroup
$\left\{  \varphi_{t}\right\}  $ on $B_{N}$ such that $\varphi_{1}=\varphi$ if
and only if there is a matrix $M\in\mathbb{C}^{N\times N}$ such that
$A=\exp\left(  M\right)  $, and for every $z\in B_{N},$
\[
\operatorname{Re}\left[  \left\langle Mz,z\right\rangle -\delta\left\langle
Mz,e_{1}\right\rangle \left\vert z\right\vert ^{2}\right]  \leq0\text{.}%
\]

\end{theorem}

The Siegel half-plane domain of $\mathbb{C}^{N}$ is defined by
\[
\mathbb{H}^{N}=\left\{  (u_{1},u^{\prime})\in\mathbb{C}\times\mathbb{C}%
^{N-1}:\operatorname{Im}u_{1}>\left\vert u^{\prime}\right\vert ^{2}\right\}
.
\]
$\mathbb{H}^{N}$ is biholomorphic to $B_{N}$ via the Cayley transformation:
\[
\sigma(z_{1},z^{\prime})=\left(  i\dfrac{1+z_{1}}{1-z_{1}},\dfrac{iz^{\prime}%
}{1-z_{1}}\right)  .
\]
Since $\sigma$ is linear fractional, for every $\varphi\in\operatorname{LFM}%
\left(  B_{N}\right)  $, there is $\psi\in\operatorname{LFM}\left(
\mathbb{H}^{N}\right)  $ such that $\left(  \varphi,B_{N}\right)  $ is
conjugated to $\left(  \psi,\mathbb{H}^{N}\right)  $.

The embedding problem of non-elliptic cases are much more complicated than the
cases of elliptic ones. Some known conclusions can be found in \cite{BCM1}. In
this article, we give two positive results first, see Theorem \ref{paraim} and
Theorem \ref{hypernormal}. Although they seems cannot be verified easily,
still there are really useful when the Jacobian of the map is normal. And then
we get the following theorems.

\begin{theorem}
\label{paramain}Let $\varphi\in\operatorname{LFM}\left(  B_{N}\right)  $ and
$\left(  \varphi,B_{N}\right)  $ is conjugated to $\left(  \psi,\mathbb{H}%
^{N}\right)  $ with%
\[
\psi\left(  z,u,v,w\right)  =\left(  z+2i\left\langle u,a\right\rangle
+2i\left\langle w,c\right\rangle +b,u+a,Dv,Aw\right)  ,
\]
where $D$ is a diagonal matrix with $\sigma\left(  D\right)  \in
\partial\mathbb{D}\backslash\left\{  1\right\}  $, $A=diag\left(  \lambda
_{1},\cdots,\lambda_{r}\right)  $ and $0<\left\vert \lambda_{j}\right\vert <1$
for $j=1,2,\cdots,r$. Let $\lambda_{j}=\exp\left(  -u_{j}+iv_{j}\right)  $
where $u_{j}>0,v_{j}\in\lbrack0,2\pi)$ for $j=1,2,\cdots,r,$ and%
\[
\Theta=diag\left(  \frac{1}{2u_{1}}\frac{\left(  u_{1}^{2}+v_{1}^{2}\right)
}{\left\vert 1-\lambda_{1}\right\vert ^{2}},\cdots,\frac{1}{2u_{r}}%
\frac{\left(  u_{r}^{2}+v_{r}^{2}\right)  }{\left\vert 1-\lambda
_{r}\right\vert ^{2}}\right)  .
\]
If%
\[
\operatorname{Im}b-\left\vert a\right\vert ^{2}\geq c^{H}\Theta c,
\]
then $\varphi$ can be embedded into a semigroup on $B_{N}$.
\end{theorem}

\begin{theorem}
\label{hypermain}Let $\varphi\in\operatorname{LFM}\left(  B_{N}\right)  $ be
hyperbolic and $\left(  \varphi,B_{N}\right)  $ is conjugated to $\left(
\psi,\mathbb{H}^{N}\right)  $ with%
\[
\psi\left(  z,u,v,w\right)  =\left(  \lambda z+2i\left\langle w,a\right\rangle
+b,\sqrt{\lambda}u,\sqrt{\lambda}Dv,\sqrt{\lambda}Aw\right)  ,
\]
where $D$ is a diagonal matrix with $\sigma\left(  D\right)  \in
\partial\mathcal{D}\backslash\left\{  1\right\}  $, $A=diag\left(  \lambda
_{1},\cdots,\lambda_{r}\right)  $ and $0<\left\vert \lambda_{j}\right\vert <1$
for $j=1,2,\cdots,r$. Let
\[
\lambda_{j}=\exp\left(  -u_{j}+iv_{j}\right)
\]
with $u_{j}>0$, $v_{j}\in\lbrack0,2\pi)$ and%
\[
\Theta=diag\left(  \frac{\lambda-1}{2u_{1}\ln\lambda}\frac{\left(  \frac
{\ln\lambda}{2}+u_{1}\right)  ^{2}+v_{1}^{2}}{\left\vert \lambda-\sqrt
{\lambda}\lambda_{1}\right\vert ^{2}},\cdots,\frac{\lambda-1}{2u_{r}\ln
\lambda}\frac{\left(  \frac{\ln\lambda}{2}+u_{r}\right)  ^{2}+v_{r}^{2}%
}{\left\vert \lambda-\sqrt{\lambda}\lambda_{r}\right\vert ^{2}}\right)  .
\]
If
\[
\operatorname{Im}b\geq\left\langle \Theta a,a\right\rangle ,
\]
then $\varphi$ can be embedded into a semigroup of $B_{N}$.
\end{theorem}

According to these Theorems, we find that an automorphism of $B_{N}$ can
always be embedded into a semigroup.

\begin{corollary}
\label{autmain}Let $\varphi$ be an automorphism of $B_{N}$, then there is a
semigroup $\left\{  \varphi_{t}\right\}  $ of automorphisms of $B_{N}$ such
that $\varphi_{1}=\varphi$.
\end{corollary}

Finally, we apply the above theorem to the case when the dimension $N=2$. The
results are very simple.

\begin{theorem}
\label{main3}Let $\varphi\in\operatorname{LFM}\left(  B_{2}\right)  $ be
parabolic. Then

\begin{enumerate}
\item[$\left(  1\right)  $] $\left(  \varphi,B_{2}\right)  $ is conjugated to
$\left(  \psi_{1},\mathbb{H}^{2}\right)  ,$ or $\left(  \psi_{2}%
,\mathbb{H}^{2}\right)  $ or $\left(  \psi_{3},\mathbb{H}^{2}\right)  $,
where
\begin{align*}
\psi_{1}\left(  u_{1},u_{2}\right)   &  =\left(  u_{1}+2ibu_{2}+c,\lambda
u_{2}\right)  ,\\
\psi_{2}\left(  u_{1,}u_{2}\right)   &  =\left(  u_{1}+c,e^{i\theta}%
u_{2}\right)  ,\\
\psi_{3}\left(  u_{1},u_{2}\right)   &  =\left(  u_{1}+2i\bar{a}u_{2}%
+c,u_{2}+a\right)  .
\end{align*}
with some specific $a,b,c\in\mathbb{C}$ and $\lambda\in\left(  0,1\right)  $.

\item[$\left(  2\right)  $] Let%
\[
\lambda=\exp\left(  -\mu+iv\right)
\]
where $\mu>0$ and $v\in\lbrack0,2\pi)$. If%
\[
\operatorname{Im}c\geq\frac{\left\vert b\right\vert ^{2}\left(  \mu^{2}%
+v^{2}\right)  }{\mu\left\vert 1-\lambda\right\vert ^{2}},
\]
then $\psi_{1}$ can be embedding into a semigroup of $\mathbb{H}^{2}$.
$\psi_{2}$ and $\psi_{3}$ can always be embedded into a semigroup on
$\mathbb{H}^{2}$.
\end{enumerate}
\end{theorem}

\begin{theorem}
\label{main4}Let $\varphi\in\operatorname{LFM}\left(  B_{2}\right)  $ be
hyperbolic, Then

\begin{enumerate}
\item[$\left(  1\right)  $] $\left(  \varphi,\mathbb{B}_{2}\right)  $ is
conjugated to $\left(  \psi_{1},\mathbb{H}^{2}\right)  $ or $\left(  \psi
_{2},\mathbb{H}^{2}\right)  $, where
\begin{align*}
\psi_{1}\left(  u_{1},u_{2}\right)   &  =\left(  \lambda u_{1}+2i\left\langle
u_{2},b\right\rangle +c,\sqrt{\lambda}\alpha u_{2}\right)  ,\\
\psi_{2}\left(  u_{1},u_{2}\right)   &  =\left(  \lambda u_{1}+a,u_{2}%
+b\right)  .
\end{align*}

\item[$\left(  2\right)  $] Let $\alpha=e^{\beta+i\gamma}$. If%
\[
\operatorname{Im}c\geq\frac{\lambda-1}{2\beta\ln\lambda}\frac{\left(
\frac{\ln\lambda}{2}+\beta\right)  ^{2}+\gamma^{2}}{\left\vert \lambda
-\sqrt{\lambda}\alpha\right\vert ^{2}}\left\vert b\right\vert ^{2},
\]
then $\psi_{1}$ can be embedded into a semigroup on $\mathbb{H}^{2}$

\item[$\left(  3\right)  $] If%
\[
\operatorname{Im}a\geq\frac{\left(  \lambda-1\right)  }{\ln^{2}\lambda
}\left\vert b\right\vert ^{2},
\]
then $\psi_{2}$ can be embedded into a semigroup on $\mathbb{H}^{2}$.
\end{enumerate}
\end{theorem}

\section{Background materials}

The following lemma is a classical result about exponential of a matrix (e.g.
see \cite[P241]{Ba}).

\begin{lemma}
\label{matrix}Given any invertible matrix $A\in\mathbb{C}^{N\times N}$, there
exists a matrix $M$ such that $\exp\left(  2\pi iM\right)  =A$. If $A$ is
triangularly blocked of some type, then so are the matrices $M$. The
eigenvalues of any two such $M$ can only differ by integers, and there is a
unique matrix $M$ whose eigenvalues have real parts in the half-open interval
$[0,1)$.
\end{lemma}

A matrix $M$ is dissipative or not can be verified by the spectra norm of
$\exp\left(  tM\right)  $ as following.

\begin{proposition}
[Phillips-Lumer, \cite{BCM1}]Let $M\in\mathbb{C}^{N\times N}$, then
$\left\Vert \exp\left(  tM\right)  \right\Vert \leq1$ for any $t\geq0$ if and
only if $M$ is dissipative.
\end{proposition}

The following proposition characterizes dissipative normal matrices. The proof
is elementary and is left to the readers as an exercise.

\begin{proposition}
\label{dissi} Let $M\in\mathbb{C}^{N\times N}$ be normal. Then $M$ is
dissipative if and only if $\left\Vert \exp\left(  M\right)  \right\Vert
\leq1$.
\end{proposition}

The following are some definitions given in \cite{BCM1}. $S\subset B_{N}$ is
called a slice of $B_{N}$ if there exists an one dimensional affine subset $V$
of $\mathbb{C}^{N}$ such that%
\[
S=B_{N}\cap V.
\]
The direction subspace $V_{S}$ of $S$ is defined by%
\[
V_{S}:=\mathrm{span}\{s-s^{\prime}:s,s^{\prime}\in S\}.
\]
For a collection of slices $\left\{  S_{j}:j=1,2,\cdots,p\right\}  $, if the
dimension of the subspace spanned by the corresponding direction subspaces
$\{V_{S_{1}},\cdots,V_{S_{p}}\}$ equals to $p$, then $\left\{  S_{j}%
:j=1,2,\cdots,p\right\}  $ is said to be linear independent. For any
$\varphi\in\operatorname{LFM}\left(  B_{N}\right)  $ and any slice $S$ of
$B_{N}$, $S$ is called an invariant slice of $\varphi$ if $\varphi\left(
S\right)  \subset S$. Let%

\[
\#\text{inv}\left(  \varphi\right)  =\dim\left(  \text{span}\left\{
V_{S}:S\text{ is an invariant slice of }\varphi\right\}  \right)  .
\]

\begin{definition}
Let $\varphi$ be an elliptic self-map of $B_{N}$, $z_{0}\in B_{N}$ be a fixed
point of $\varphi$. $L_{U}\left(  \varphi,z_{0}\right)  \subset\mathbb{C}^{N}$
is called the unitary space of $\varphi$ at $z_{0}$ if%
\[
L_{U}\left(  \varphi,z_{0}\right)  =\underset{\left\vert \lambda\right\vert
=1}{%
{\displaystyle\bigoplus}
}\ker\left(  d\varphi_{z_{0}}-\lambda E\right)  ^{N}.
\]
And $u\left(  \varphi,z_{0}\right)  =\dim L_{U}\left(  \varphi,z_{0}\right)  $
is said to be the unitary index of $\varphi$ at $z_{0}.$
\end{definition}

According to Lemma 3.1 in \cite{BCM1}, $u\left(  \varphi,z_{0}\right)  \ $is
independent of the choice of the fixed point $z_{0}$. Thus the unitary index
of $\varphi$ can be denoted by $u\left(  \varphi\right)  $.

The following example shows that there exist non-linear fractional semigroups
with non-trivial linear iterates.

\begin{example}
\label{counterexample}For $k\geq\frac{16}{9}\pi$, let%
\[
\varphi_{t}\left(  z,w\right)  =\left[
\begin{array}
[c]{c}%
\exp\left(  2\pi it\right)  z\\
\frac{4+\exp\left(  4\pi it\right)  z^{2}}{\exp\left(  kt\right)  \left(
4+z^{2}\right)  }w
\end{array}
\right]  ,
\]
then $\varphi_{1}\left(  z,w\right)  =\left(  z,e^{-k}w\right)  $ is linear,
and $\varphi_{t}\left(  z,w\right)  $ is a semigroup on $\mathbb{B}_{2}.$
\end{example}

\begin{proof}
We only need to show that for any $\left(  z,w\right)  \in\mathbb{B}_{2}$,
$\varphi_{t}\left(  z,w\right)  \in\mathbb{B}_{2}$. Let%
\[
f\left(  t\right)  =\left\vert 4+\exp\left(  4\pi it\right)  z^{2}\right\vert
^{2},g\left(  t\right)  =\left\vert 4+z^{2}\right\vert ^{2}\exp\left(
2kt\right)  ,
\]
then $f\left(  0\right)  =g\left(  0\right)  =\left\vert 4+z^{2}\right\vert
^{2}$. And%
\begin{align*}
f^{\prime}\left(  t\right)   &  =2\operatorname{Re}\left(  16\pi iz^{2}%
\exp\left(  4\pi it\right)  \right)  \leq32\pi,\\
g^{\prime}\left(  t\right)   &  =\left\vert 4+z^{2}\right\vert ^{2}\exp\left(
2kt\right)  \cdot2k\geq18k.
\end{align*}
Therefore, $f\left(  t\right)  \leq g\left(  t\right)  $, and as an
conclusion, for any $\left(  z,w\right)  \in\mathbb{B}_{2}$, $\varphi
_{t}\left(  z,w\right)  \in\mathbb{B}_{2}$.
\end{proof}

The following results show that there exists some semigroups have the property
that if one iterate is linear, then all the other iterates are linear.

\begin{lemma}
\label{block}Let $d_{0},\cdots,d_{m}\in\mathbb{N}$, $F:\mathbb{C}^{d_{1}%
}\times\cdots\times\mathbb{C}^{d_{m}}\rightarrow\mathbb{C}^{d_{0}}$ be a
multilinear mapping, $\lambda_{0},\lambda_{1},\lambda_{2},\cdots,\lambda
_{m}\in\mathbb{\mathbb{C}}\backslash\left\{  0\right\}  $ with $\left\vert
\lambda_{0}\lambda_{1}^{-1}\cdots\lambda_{m}^{-1}\right\vert \geq1$. Let%
\[
J_{d_{i}}\left(  \lambda_{i}\right)  =\left[
\begin{array}
[c]{cccc}%
\lambda_{i} & 1 &  & \\
& \ddots & \ddots & \\
&  & \lambda_{i} & 1\\
&  &  & \lambda_{i}%
\end{array}
\right]  _{d_{i}\times d_{i}}\in\mathbb{C}^{d_{i}\times d_{i}},i=0,\cdots,m.
\]
If for every $\left(  v_{1},\cdots,v_{m}\right)  \in\mathbb{C}^{d_{1}}%
\times\cdots\times\mathbb{C}^{d_{m}},$
\[
J_{d_{0}}\left(  \lambda_{0}\right)  ^{n}\circ F\circ\left(  J_{d_{1}}\left(
\lambda_{1}\right)  ^{-n}v_{1},\cdots,J_{d_{m}}\left(  \lambda_{m}\right)
^{-n}v_{m}\right)  \rightarrow0\text{, }\left(  n\rightarrow+\infty\right)  ,
\]
then $F\equiv0.$
\end{lemma}

\begin{proof}
First of all,
\begin{align*}
J_{d_{0}}\left(  \lambda_{0}\right)  ^{n}  &  =\left(  J_{d_{0}}\left(
\lambda_{0}\right)  -\lambda_{0}I+\lambda_{0}I\right)  ^{n}\\
&  =\sum_{k=0}^{d_{0}-1}C_{n}^{k}\lambda_{0}^{n-k}\left(  J_{d_{0}}\left(
\lambda_{0}\right)  -\lambda_{0}I\right)  ^{k},
\end{align*}
and%
\begin{align*}
J_{d_{i}}\left(  \lambda_{i}\right)  ^{-n}  &  =\left(  J_{d_{i}}\left(
\lambda_{i}\right)  -\lambda_{i}I+\lambda_{i}I\right)  ^{-n}\\
&  =\lambda_{i}^{-n}\left(  I+\frac{J_{d_{i}}\left(  \lambda_{i}\right)
-\lambda_{i}I}{\lambda_{i}}\right)  ^{-n}\\
&  =\lambda_{i}^{-n}\sum_{k=0}^{+\infty}\frac{\left(  -n\right)  \left(
-n-1\right)  \cdots\left(  -n-k+1\right)  }{k!}\left(  \frac{J_{d_{i}}\left(
\lambda_{i}\right)  -\lambda_{i}I}{\lambda_{i}}\right)  ^{k}\\
&  =\sum_{k=0}^{d_{i}-1}\left(  -1\right)  ^{k}\frac{n\cdots\left(
n+k-1\right)  }{k!}\left(  J_{d_{i}}\left(  \lambda_{i}\right)  -\lambda
_{i}\right)  ^{k}\lambda_{i}^{-n-k}.
\end{align*}

Let $c_{k}^{\left(  n\right)  }=C_{n}^{k-1}\lambda_{0}^{n-k+1},d_{k,i}%
^{\left(  n\right)  }=\left(  -1\right)  ^{k-1}\frac{n\cdots\left(
n+k-2\right)  }{k!}\lambda_{i}^{-n-k+1}.$ Then%
\[
c_{k}^{\left(  n\right)  }\sim n^{k-1}\lambda_{0}^{n},\;d_{k,i}^{\left(
n\right)  }\sim n^{k-1}\lambda_{i}^{-n}.
\]
For $k=1,\cdots,m$, let $\left\{  e_{1}^{\left(  d_{k}\right)  }%
,\cdots,e_{d_{k}}^{\left(  d_{k}\right)  }\right\}  $ be the standard basis of
$\mathbb{C}^{d_{k}}$. Since $F$ is a multilinear mapping, we may assume that
$F=\left(  F_{1},\cdots,F_{d_{0}}\right)  ^{T}$, and that%
\[
F_{d_{0}}=\sum_{i_{1},\cdots,i_{m}}a^{i_{1},\cdots,i_{m}}e_{i_{1}}^{\left(
d_{1}\right)  }\otimes\cdots\otimes e_{i_{m}}^{\left(  d_{m}\right)  },
\]
where $e_{i_{1}}^{\left(  d_{1}\right)  }\otimes\cdots\otimes e_{i_{m}%
}^{\left(  d_{m}\right)  }$ is tensor product of $e_{i_{1}}^{\left(
d_{1}\right)  },\cdots,e_{i_{m}}^{\left(  d_{m}\right)  }.$ By assumption,%
\[
\lambda_{0}^{n}\sum_{i_{1},\cdots,i_{m}}a^{i_{1},\cdots,i_{m}}\left(
J_{d_{1}}\left(  \lambda_{1}\right)  ^{-n}e_{i_{1}}^{\left(  d_{1}\right)
}\right)  \otimes\cdots\otimes\left(  J_{d_{m}}\left(  \lambda_{m}\right)
^{-n}e_{i_{m}}^{\left(  d_{m}\right)  }\right)  \rightarrow0.
\]
Since%
\begin{align*}
&  \sum_{i_{1},\cdots,i_{m}}a^{i_{1},\cdots,i_{m}}\left(  J_{d_{1}}\left(
\lambda_{1}\right)  ^{-n}e_{i_{1}}^{\left(  d_{1}\right)  }\right)
\otimes\cdots\otimes\left(  J_{d_{m}}\left(  \lambda_{m}\right)  ^{-n}%
e_{i_{m}}^{\left(  d_{m}\right)  }\right) \\
&  =\sum_{i_{1},\cdots,i_{m}}a^{i_{1},\cdots,i_{m}}\left(  \sum_{j_{1}%
=1}^{i_{1}}d_{j_{1},1}^{\left(  n\right)  }e_{j_{1}}^{\left(  d_{1}\right)
}\right)  \otimes\cdots\otimes\left(  \sum_{j_{m}=1}^{i_{m}}d_{j_{m}%
,m}^{\left(  n\right)  }e_{j_{M}}^{\left(  d_{m}\right)  }\right) \\
&  =\left(  \sum_{j_{1}=1}^{d_{1}}\cdots\sum_{j_{m}=1}^{d_{m}}\right)  \left(
\sum_{i_{1}=j_{1}}^{d_{1}}\cdots\sum_{i_{m}=j_{m}}^{d_{m}}a^{i_{1}%
,\cdots,i_{m}}d_{j_{1},1}^{\left(  n\right)  }\cdots d_{j_{m},m}^{\left(
n\right)  }e_{j_{1}}^{\left(  d_{1}\right)  }\otimes\cdots\otimes e_{j_{M}%
}^{\left(  d_{m}\right)  }\right)  ,
\end{align*}
we have for any $\left(  j_{1},\cdots,j_{m}\right)  ,$%
\[
\lambda_{0}^{n}\sum_{i_{1}=j_{1}}^{d_{1}}\cdots\sum_{i_{m}=j_{m}}^{d_{m}%
}a^{i_{1},\cdots,i_{m}}d_{j_{1},1}^{\left(  n\right)  }\cdots d_{j_{m}%
,m}^{\left(  n\right)  }\rightarrow0\text{.}%
\]
When $j_{1}=d_{1},\cdots,j_{m}=d_{m},$ the above formula implies that%
\[
\left\vert \lambda_{0}^{n}a^{d_{1},\cdots,d_{m}}d_{d_{1},1}^{\left(  n\right)
}\cdots d_{d_{m},m}^{\left(  n\right)  }\right\vert \sim\left\vert
a^{d_{1},\cdots,d_{m}}\right\vert \cdot n^{d_{1}+\cdots+d_{m}-m}\cdot\left(
\lambda_{0}\lambda_{1}^{-1}\cdots\lambda_{m}^{-1}\right)  ^{n}\rightarrow
0\text{.}%
\]
Since $n^{d_{1}+\cdots+d_{m}-m}\cdot\left(  \lambda_{0}\lambda_{1}^{-1}%
\cdots\lambda_{m}^{-1}\right)  ^{n}\rightarrow+\infty$, we find that
$a^{d_{1},\cdots,d_{m}}=0.$ When $j_{1}=d_{1},\cdots,j_{m}=d_{m}-1$,
\begin{align*}
&  \lambda_{0}^{n}\left(  a^{d_{1},\cdots,d_{m}}d_{d_{1},1}^{\left(  n\right)
}\cdots d_{d_{m},m}^{\left(  n\right)  }+a^{d_{1},\cdots,d_{m}-1}d_{d_{1}%
,1}^{\left(  n\right)  }\cdots d_{d_{m}-1,m}^{\left(  n\right)  }\right) \\
&  =\lambda_{0}^{n}a^{d_{1},\cdots,d_{m}-1}d_{d_{1},1}^{\left(  n\right)
}\cdots d_{d_{m}-1,m}^{\left(  n\right)  }\rightarrow0,
\end{align*}
and this formula implies that $a^{d_{1},\cdots,d_{m}-1}=0$. By induction,
$a^{j_{1},\cdots,j_{m}}=0$ for all $\left(  j_{1},\cdots,j_{m}\right)  $, and
$F_{d_{0}}\equiv0$.

Notice that
\[
J_{d_{0}}\left(  \lambda_{0}\right)  ^{n}\circ F=\left(  \cdots,\lambda
_{0}^{n}F_{d_{0}-1},0\right)  ^{T},
\]
using the same method as above, $F_{d_{0}-1}\equiv0$, then $F_{d_{0}-2}%
\equiv0$, and finally, $F_{1}\equiv0$. As a conclusion, $F\equiv0$.
\end{proof}

\begin{lemma}
\label{lemmamain}Let $D\subset\mathbb{C}^{N}$ be an open domain with $0\in
\bar{D}$.
\[
J=\operatorname{diag}\left(  J_{d_{1}}\left(  \lambda_{1}\right)
,\cdots,J_{d_{m}}\left(  \lambda_{m}\right)  \right)  \in\mathbb{C}^{N\times
N}%
\]
with $1>\left\vert \lambda_{1}\right\vert \geq\left\vert \lambda
_{2}\right\vert \cdots\geq\left\vert \lambda_{m}\right\vert >0,$ and $J\left(
D\right)  \subset D$. Denote
\[
K=\min\left\{  k\in\mathbb{N}:\left\vert \lambda_{m}\lambda_{1}^{-k}%
\right\vert \geq1\right\}  .
\]
Let $\varphi:D\rightarrow D$ be holomorphic and $\varphi\circ J=J\circ\varphi
$.$\ $If $\varphi$ is of class $C^{K}$ at $0$, then every component of
$\varphi$ is a polynomial with degree no more than $K+1$.
\end{lemma}

\begin{proof}
Since $\varphi\circ J=J\circ\varphi$, for any $v_{1},\cdots,v_{K}\in
\mathbb{C}^{N}$, and any $z\in D$,%
\[
\varphi^{\left(  K\right)  }\left(  Jz\right)  \left(  Jv_{1},\cdots
,Jv_{K}\right)  =J\circ\varphi^{\left(  K\right)  }\left(  z\right)  \left(
v_{1},\cdots,v_{K}\right)  \text{.}%
\]
Therefore,%
\[
\varphi^{\left(  K\right)  }\left(  Jz\right)  \left(  v_{1},\cdots
,v_{K}\right)  =J\circ\varphi^{\left(  K\right)  }\left(  z\right)  \left(
J^{-1}v_{1},\cdots,J^{-1}v_{K}\right)  ;
\]
let $z=0,$ $\varphi^{\left(  K\right)  }\left(  0\right)  \left(  v_{1}%
,\cdots,v_{k}\right)  =J\circ\varphi^{\left(  K\right)  }\left(  0\right)
\left(  J^{-1}v_{1},\cdots,J^{-1}v_{K}\right)  $. Now%
\begin{align*}
J^{n}\circ\varphi^{\left(  K\right)  }\left(  z\right)  \left(  J^{-n}%
v_{1},\cdots,J^{-n}v_{K}\right)   &  =\varphi^{\left(  K\right)  }\left(
J^{n}z\right)  \left(  v_{1},\cdots,v_{K}\right) \\
&  \rightarrow\varphi^{\left(  K\right)  }\left(  0\right)  \left(
v_{1},\cdots,v_{K}\right)  \text{,}%
\end{align*}
we have%
\[
J^{n}\circ\left(  \varphi^{\left(  K\right)  }\left(  z\right)  -\varphi
^{\left(  K\right)  }\left(  0\right)  \right)  \left(  J^{-n}v_{1}%
,\cdots,J^{-n}v_{K}\right)  \rightarrow0\text{, }\left(  n\rightarrow
+\infty\right)  \text{.}%
\]
Let $F=\varphi^{\left(  K\right)  }\left(  z\right)  -\varphi^{\left(
K\right)  }\left(  0\right)  $, then $F:\mathbb{C}^{N}\times\cdots
\times\mathbb{C}^{N}\rightarrow\mathbb{C}^{N}$ is a multilinear mapping.
Divide $F$ into $m^{K+1}$ blocks: $F_{j_{1},\cdots,j_{K}}^{i}:i,j_{1}%
,\cdots,j_{K}=1\cdots,m$, where%
\[
F_{j_{1},\cdots,j_{K}}^{i}:\mathbb{C}^{d_{j_{1}}}\times\cdots\mathbb{C}%
^{d_{j_{K}}}\rightarrow\mathbb{C}^{d_{i}}\text{.}%
\]
Then for any $\left(  w_{1},\cdots,w_{K}\right)  \in\mathbb{C}^{d_{j_{1}}%
}\times\cdots\mathbb{C}^{d_{j_{K}}},$%
\[
J_{d_{i}}\left(  \lambda_{i}\right)  ^{n}F_{j_{1},\cdots,j_{K}}^{i}\left(
J_{d_{j_{1}}}\left(  \lambda_{j_{1}}\right)  ^{-n}w_{1},\cdots,J_{d_{j_{K}}%
}\left(  \lambda_{j_{K}}\right)  ^{-n}w_{K}\right)  \rightarrow0,\left(
n\rightarrow0\right)  \text{.}%
\]
Since $\left\vert \lambda_{i}\lambda_{j_{1}}^{-1}\cdots\lambda_{j_{K}}%
^{-1}\right\vert \geq\left\vert \lambda_{m}\lambda_{1}^{-K}\right\vert \geq1$,
according to Lemma \ref{block}, $F_{j_{1},\cdots,j_{K}}^{i}\equiv0$, and
consequently, $F\equiv0$. Therefore, for any $z\in\mathbb{D}$,%
\[
\varphi^{\left(  K\right)  }\left(  z\right)  =\varphi^{\left(  K\right)
}\left(  0\right)  \text{,}%
\]
and as a conclusion, every component of $\varphi$ is a polynomial with degree
no more than $K+1$.
\end{proof}

\begin{theorem}
\label{linear}Let $D\subset\mathbb{C}^{N}$ be an open domain with
$0\in\overline{D},$ $A\in\mathbb{C}^{N\times N}$ is invertible with $A\left(
D\right)  \subset D,$ and $\rho\left(  A\right)  <1$. Let
\[
K=\min\left\{  k\in\mathbb{N}:\rho\left(  A^{-1}\right)  \rho\left(  A\right)
^{k}\geq1\right\}  .
\]
Suppose that $\left(  \varphi_{t}\right)  $ is a semigroup on $D$,
$\varphi_{1}\left(  z\right)  =Az,$ and for every $t>0$, $\varphi_{t}$ is of
class $C^{K}$ at $0$. Then for every $t>0$, $\varphi_{t}$ is a linear map.
\end{theorem}

\begin{proof}
There is a invertible matrix $P$ such that $P^{-1}AP$ is a Jordan matrix. Let
$\psi_{t}=P^{-1}\circ\varphi_{t}\circ P$. Then $\left\{  \psi_{t}\right\}  $
is a semigroup on $P^{-1}\left(  D\right)  $. According to Lemma
\ref{lemmamain}, for every $t>0$, every component of $\psi_{t}$ is a
polynomial with degree not more than $K+1$. If $\psi_{t}$ is not linear, then
$\deg\psi_{t}^{n}\rightarrow+\infty\left(  n\rightarrow\infty\right)  $, which
is impossible. As a conclusion, for any $t>0$, $\varphi_{t}$ is linear.
\end{proof}

\section{The elliptic cases}

Let $\varphi\in\operatorname{LFM}\left(  B_{N}\right)  $ be elliptic and
$z_{0}\in B_{N}$ be one of the fixed points of $\varphi$. Let $\varphi_{z_{0}%
}$ be the automorphism of $B_{N}$ such that $\varphi_{z_{0}}\left(  0\right)
=z_{0},$ $\varphi_{z_{0}}^{-1}=\varphi_{z_{0}}$. Let $\psi=\varphi_{z_{0}%
}\circ\varphi\circ\varphi_{z_{0}}$, then $\psi\in\operatorname{LFM}\left(
B_{N}\right)  $, and $\psi\left(  0\right)  =0$. There exist some
$A\in\mathbb{C}^{N\times N}$, $C\in\mathbb{C}^{N}$ such that
\begin{equation}
\psi\left(  z\right)  =\frac{Az}{\left\langle z,C\right\rangle +1}.
\label{elf}%
\end{equation}
It is easy to see that
\[
\psi^{\prime}\left(  0\right)  =\varphi_{z_{0}}^{\prime}\left(  0\right)
^{-1}\circ\varphi^{\prime}\left(  z_{0}\right)  \circ\varphi_{z_{0}}^{\prime
}\left(  0\right)  .
\]
Let $F=Fix\left(  \varphi\right)  $ be the collection of fixed points of
$\varphi$. According to \cite{R}, $F$ is the intersection of $B_{N}$ and some
subspace of $\mathbb{C}^{N}$. Let $p$ denote the dimension of this subspace
and let $u=u\left(  \varphi\right)  $ be the unitary index of $\varphi$. It is
clear that if $\left(  \varphi,B_{N}\right)  $ is conjugated to $\left(
\psi,B_{N}\right)  $, where $\psi$ is given in $\left(  \ref{elf}\right)  $,
then $u\left(  \varphi\right)  >0$ if and only if $\rho\left(  A\right)  =1$.

\begin{proposition}
$\label{algebramul}$Let $A\in\mathbb{C}^{N\times N}$, and $\left\Vert
A\right\Vert \leq1$. If $\lambda$ is an eigenvalue of $A$ with $\left\vert
\lambda\right\vert =1$, then the generalized eigenspace of $\lambda$ coincides
with the eigenspace of $\lambda$.
\end{proposition}

\begin{proof}
According to Schur's Triangularization Theorem (See, for instance,
\cite[P508]{Me}), there is a unitary matrix $U\in\mathbb{C}^{N\times N}$ and
an upper-triangular matrix%
\[
T=\left[
\begin{array}
[c]{cccc}%
\lambda & a_{12} & \cdots & a_{1N}\\
0 & \lambda_{2} & \ast & a_{2N}\\
\vdots & \vdots & \ddots & \vdots\\
0 & 0 & \cdots & \lambda_{N}%
\end{array}
\right]  ,
\]
such that%
\[
A=U^{H}TU.
\]
Then $\left\Vert T\right\Vert \leq1$ since $\left\Vert A\right\Vert \leq1$.

Suppose there is a subscript $i$ such that $a_{1i}\neq0$. Let
\[
z_{\lambda}=\left(  \frac{\bar{\lambda}}{\sqrt{1+\left\vert a_{1i}\right\vert
^{2}}},0,\cdots,\frac{\bar{a}_{1i}}{\sqrt{1+\left\vert a_{1i}\right\vert ^{2}%
}},\cdots,0\right)  ^{T}\in\mathbb{C}^{N}.
\]
Then $z_{\lambda}$ is a unit vector and
\[
Tz_{\lambda}=\left(  \frac{\bar{\lambda}\lambda}{\sqrt{1+\left\vert
a_{1i}\right\vert ^{2}}}+\frac{\bar{a}_{1i}a_{1i}}{\sqrt{1+\left\vert
a_{1i}\right\vert ^{2}}},\cdots\right)  ^{T}=\left(  \sqrt{1+\left\vert
a_{1i}\right\vert ^{2}},\cdots\right)  .
\]
Therefore
\[
\left\vert Tz_{\lambda}\right\vert \geq\sqrt{1+\left\vert a_{1i}\right\vert
^{2}},
\]
which is impossible since $\left\Vert T\right\Vert \leq1$. As a consequence,
$a_{1i}=0$ for all $i=2,\cdots,n.$ This completes the proof of this lemma.
\end{proof}

The following theorem characterize elliptic linear fractional maps on $B_{N}$.

\begin{theorem}
\label{ell_normal}Let $\varphi\in\operatorname{LFM}\left(  B_{N}\right)  $ be elliptic.

\begin{enumerate}
\item If $u=u\left(  \varphi\right)  >0$, then $\left(  \varphi,B_{N}\right)
$ is conjugated to $\left(  \psi,B_{N}\right)  $ defined by%
\[
\psi\left(  z^{\prime},z^{\prime\prime}\right)  =\left(  \Lambda z^{\prime
},A_{1}z^{\prime\prime}\right)  ,
\]
where $\left(  z^{\prime},z^{\prime\prime}\right)  \in\mathbb{C}^{u}%
\times\mathbb{C}^{N-u}\cap B_{N}$, $\Lambda$ is a diagonal and unitary matrix
of order $u$ and $A_{1}$ is a matrix of order $N-u$ with $\rho\left(
A_{1}\right)  <1$, $\left\Vert A_{1}\right\Vert \leq1$.

\item If $u=u\left(  \varphi\right)  =0$. Then $\left(  \varphi,B_{N}\right)
$ is conjugated to $\left(  \psi,B_{N}\right)  $ defined by%
\[
\psi\left(  z\right)  =\frac{Az}{\delta\left\langle z,\left(  A^{H}-E\right)
e_{1}\right\rangle +1},
\]
where $A$ is a matrix of order $N$ with $\rho\left(  A\right)  <1$,
$\left\Vert A\right\Vert \leq1$ and $\delta\in\left[  0,1\right]  $,
$e_{1}=\left(  1,0,\cdots,0\right)  ^{T}$. Moreover, there is domain $D$ with
$0\in D$ which is biholomorphic equivalent to $B_{N}$ such that $\left(
\varphi,B_{N}\right)  $ is conjugated to $\left(  \tilde{\psi},D\right)  $
with%
\[
\tilde{\psi}\left(  z\right)  =Az.
\]

\end{enumerate}
\end{theorem}

\begin{proof}
\begin{enumerate}
\item Resulting from the previous discussion, we may assume that%
\[
\varphi\left(  z\right)  =\frac{Az}{\left\langle z,C\right\rangle +1}.
\]
Simple computation indicates that the Jacobi matrix of $\varphi$ at the origin
$d\varphi_{O}=A$. According to Schwartz's lemma (see \cite{R}), we have
$\left\Vert A\right\Vert \leq1$. Due to Proposition \ref{algebramul} and
$\rho\left(  A\right)  =1$, there is a unitary matrix $U$, such that%
\[
U^{H}AU=\left[
\begin{array}
[c]{cc}%
\Lambda & \\
& A_{1}%
\end{array}
\right]  ,
\]
where%
\[
\Lambda=\left[
\begin{array}
[c]{cccc}%
\lambda_{1} &  &  & \\
& \lambda_{2} &  & \\
&  & \ddots & \\
&  &  & \lambda_{u}%
\end{array}
\right]  ,
\]
and $\left\vert \lambda_{j}\right\vert =1$ for $j=1,2,\cdots,u$. $A_{1}$ is a
matrix of order $N-u$ with $\left\Vert A_{1}\right\Vert \leq1$ and
$\rho\left(  A_{1}\right)  <1$.

Let $\psi\left(  z\right)  =U^{H}\left(  \varphi\left(  Uz\right)  \right)  $.
Then%
\[
\psi\left(  z^{\prime},z^{\prime\prime}\right)  =\frac{U^{H}AUz}{\left\langle
z,U^{H}C\right\rangle +1}=\frac{\left(  \Lambda z^{\prime},A_{1}%
z^{\prime\prime}\right)  }{\left\langle z,U^{H}C\right\rangle +1},
\]
where $\left(  z^{\prime},z^{\prime\prime}\right)  \in\mathbb{C}^{u}%
\times\mathbb{C}^{N-u}\cap B_{N}$. We denote%
\[
U^{H}C=\left(  c^{\prime},c^{\prime\prime}\right)  \in\mathbb{C}^{u}%
\times\mathbb{C}^{N-u},
\]
thus%
\[
\psi\left(  z^{\prime},z^{\prime\prime}\right)  =\frac{\left(  \Lambda
z^{\prime},A_{1}z^{\prime\prime}\right)  }{\left\langle z^{\prime},c^{\prime
}\right\rangle +\left\langle z^{\prime\prime},c^{\prime\prime}\right\rangle
+1}.
\]
Since $\psi\left(  B_{N}\right)  \subset B_{N}$, we find that
\[
\psi\left(  \left\{  \left(  z^{\prime},0\right)  \in\mathbb{\mathbb{C}}%
^{N}:\left\vert z^{\prime}\right\vert <1\right\}  \right)  \subset\left\{
\left(  z^{\prime},0\right)  \in\mathbb{C}^{N}:\left\vert z^{\prime
}\right\vert <1\right\}  ,
\]
and consequently $\psi_{1}\left(  z^{\prime}\right)  \overset{\Delta}{=}%
\psi\left(  z^{\prime},O\right)  =\frac{\Lambda z^{\prime}}{\left\langle
z^{\prime},c^{\prime}\right\rangle +1}$ is a self-map of the unit ball of
$\mathbb{C}^{u}$ and%
\[
d\left(  \psi_{1}\right)  _{O}=\Lambda.
\]
By Schwartz's lemma on the ball, $\psi_{1}$ is linear and as a consequence
$c^{\prime}=0$. As a result,%
\[
\psi\left(  z\right)  =\dfrac{\left(  \Lambda z^{\prime},A_{1}z^{\prime\prime
}\right)  }{\left\langle z^{\prime\prime},c^{\prime\prime}\right\rangle +1}.
\]
Since $\psi\in\operatorname{LFM}\left(  B_{N}\right)  $, $\left\vert
c^{\prime\prime}\right\vert <1$. If $c^{\prime\prime}\neq0,$ let
\[
z_{t}=\left(  \sqrt{1-t^{2}\left\vert c^{\prime\prime}\right\vert ^{2}}%
e_{1}^{\prime},-tc^{\prime\prime}\right)  ,
\]
where $e_{1}^{\prime}=\left(  1,0,\cdots,0\right)  ^{T}\in\mathbb{C}^{u}$.
Then for all $t\in\left[  0,1\right]  $, $\left\vert z_{t}\right\vert =1$, and%
\[
\left\vert \psi\left(  z_{t}\right)  \right\vert ^{2}=\left\vert
\dfrac{\left(  \sqrt{1-t^{2}\left\vert c^{\prime\prime}\right\vert ^{2}%
}\Lambda e_{1}^{\prime},-tA_{1}c^{\prime\prime}\right)  }{1-t\left\vert
c^{\prime\prime}\right\vert ^{2}}\right\vert ^{2}\geq\frac{1-t^{2}\left\vert
c^{\prime\prime}\right\vert ^{2}}{\left(  1-t\left\vert c^{\prime\prime
}\right\vert ^{2}\right)  ^{2}}.
\]
It is very easy to see that $\left\vert \psi\left(  z_{t}\right)  \right\vert
^{2}$ is increasing with respect to $t\in(0,1]$. Therefore for any $t\in
(0,1]$,
\[
\left\vert \psi\left(  z_{t}\right)  \right\vert ^{2}>\lim_{t\rightarrow
0}\frac{1-t^{2}\left\vert c^{\prime\prime}\right\vert ^{2}}{\left(
1-t\left\vert c^{\prime\prime}\right\vert ^{2}\right)  ^{2}}=1.
\]
It is impossible since $\psi\left(  \overline{B_{N}}\right)  \subset
\overline{B_{N}}$. As a consequence, $c^{\prime\prime}=O$ and%
\[
\psi\left(  z\right)  =\left(  \Lambda z^{\prime},A_{1}z^{\prime\prime
}\right)  .
\]

\item See the proof in \cite[Proposition 3.4]{kms2011}
\end{enumerate}
\end{proof}

We will make use of the following generalization of Berkson-Porta's criterion
due to Aharonov, Elin, Reich and Shoikhet (see Theorem 1.3, \cite{B1}):

\begin{lemma}
\label{infinity1} Let $F:B_{N}\rightarrow\mathbb{C}^{N}$ be holomorphic. $F$
is the infinitesimal generator of a semigroup of holomorphic self-maps of
$B_{N}$ fixing the origin if and only if
\[
F\left(  z\right)  =-Q\left(  z\right)  z,
\]
where $Q\left(  z\right)  $ is a matrix of order $N$ with holomorphic entries
such that%
\[
\operatorname{Re}\left\langle Q\left(  z\right)  ,z\right\rangle \geq0.
\]

\end{lemma}

Now we can prove Theorem \ref{main1} and Theorem \ref{main2}.

\begin{proof}
[Proof of Theorem \ref{main1}]It is easy to see that there is a real diagonal
matrix $\Theta$ such that $\exp\left(  i\Theta\right)  =\Lambda$.

Suppose firstly that there is a dissipative matrix $M$ such that $\exp\left(
M\right)  =A_{1}$. Let%
\[
\varphi_{t}\left(  z^{\prime},z^{\prime\prime}\right)  =\left(  \exp\left(
it\Theta\right)  z^{\prime},\exp\left(  tM\right)  z^{\prime\prime}\right)  .
\]
Then $\varphi_{t}\left(  B_{N}\right)  \subset B_{N}$ since $\left\Vert
\exp\left(  it\Theta\right)  \right\Vert =1$ and $\left\Vert \exp\left(
tM\right)  \right\Vert \leq1$. Moreover, $\varphi_{t+s}=\varphi_{t}%
\circ\varphi_{s}=\varphi_{s}\circ\varphi_{t}$. As a result, $\left\{
\varphi_{t}\right\}  $ is a semigroup of $B_{N}$ with $\varphi_{1}=\varphi.$

On the other hand, if $\varphi$ can be embedded into a semigroup of linear
fractional self-maps $\left\{  \varphi_{t}\right\}  $, $\varphi_{t}$ is
conjugated to the following linear fractional map owning to Theorem 3.2 of
\cite{BCM1}:%
\[
\psi_{t}\left(  z^{\prime},z^{\prime\prime}\right)  =\left(  \exp\left(
it\tilde{\Theta}\right)  z^{\prime},\exp\left(  t\tilde{M}\right)
z^{\prime\prime}\right)  ,
\]
where $\tilde{\Theta}$ is a real diagonal matrix, $\tilde{M}$ is dissipative
and all eigenvalues locates on the left half plane. Suppose that
$\varphi_{t_{0}}=$ $\varphi$. Then $\varphi$ is conjugated to%
\[
\psi_{t_{0}}\left(  z^{\prime},z^{\prime\prime}\right)  =\left(  \exp\left(
it_{0}\tilde{\Theta}\right)  z^{\prime},\exp\left(  t_{0}\tilde{M}\right)
z^{\prime\prime}\right)  .
\]
Since $\psi_{t_{0}}^{\prime}\left(  0\right)  $ and $\varphi_{t_{0}}^{\prime
}\left(  0\right)  $ are similar, there exist matrixes $U$ and $V$ such that%
\[
\Lambda=U^{-1}\exp\left(  it_{0}\tilde{\Theta}\right)  U\text{ and }%
A_{1}=V^{-1}\exp\left(  t_{0}\tilde{M}\right)  V=\exp\left(  t_{0}V^{-1}%
\tilde{M}V\right)  .
\]
Let $M=t_{0}V^{-1}\tilde{M}V$, then $M$ is dissipative, and $\sigma\left(
M\right)  \subset\mathbb{R}^{-}\cup\left\{  0\right\}  $.
\end{proof}

\begin{proof}
[Proof of Theorem \ref{main2}]If $\varphi$ can be embedded into $\left\{
\varphi_{t}\right\}  $ which is a semigroup on $B_{N}$, according to Theorem
\ref{ell_normal}, there is $A\in\mathbb{C}^{N\times N}$ with $\rho\left(
A\right)  <1,\left\Vert A\right\Vert \leq1$, a domain $D$ and a linear
fractional map $\tau:B_{N}\rightarrow D$ such that $\tau\circ\varphi\circ
\tau^{-1}\left(  z\right)  =Az$. If there is a semigroup on $B_{N}$ such that
$\varphi_{1}=\varphi$, then $\left\{  \tau\circ\varphi_{t}\circ\tau
^{-1}\right\}  $ is a semigroup on $D$ and $\tau\circ\varphi_{t}\circ\tau
^{-1}$ is holomorphic at $0$. According to Theorem \ref{linear}, $\tau
\circ\varphi_{t}\circ\tau^{-1}$ is linear for all $t>0$. Since $\tau$ and
$\tau^{-1}$ is linear fractional (see proof in \cite{kms2011}), we see that
$\left\{  \varphi_{t}\right\}  $ is a semigroup of linear fractional maps. Due
to Theorem 3.2 of \cite{BCM1}, there is a matrix $M$ such that
\[
\varphi_{t}\left(  z\right)  =\frac{\exp\left(  tM\right)  z}{\delta
\left\langle z,\left(  \exp\left(  tM\right)  ^{H}-I\right)  e_{1}%
\right\rangle +1}.
\]
$A=\exp\left(  M\right)  $ for $\varphi=\varphi_{1}$. Easy computation gives%
\[
\frac{d}{dt}\varphi_{t}\left(  z\right)  =\left(  M-\delta\left\langle
M\varphi_{t}\left(  z\right)  ,e_{1}\right\rangle I\right)  \varphi_{t}\left(
z\right)  .
\]
Thus
\[
F\left(  z\right)  =-\left(  M-\delta\left\langle Mz,e_{1}\right\rangle
I\right)  z
\]
is the infinitesimal generator of $\left\{  \varphi_{t}\right\}  $ . By Lemma
\ref{infinity1}, we obtain%
\[
\operatorname{Re}\left\langle \left(  M-\delta\left\langle Mz,e_{1}%
\right\rangle I\right)  z,z\right\rangle =\operatorname{Re}\left[
\left\langle Mz,z\right\rangle -\delta\left\langle Mz,e_{1}\right\rangle
\left\vert z\right\vert ^{2}\right]  \leq0.
\]

On the other hand, if there is a matrix $M$ such that $A=\exp\left(  M\right)
$, and
\[
\operatorname{Re}\left[  \left\langle Mz,z\right\rangle -\delta\left\langle
Mz,e_{1}\right\rangle \left\vert z\right\vert ^{2}\right]  \leq0\text{.}%
\]
Since $F\left(  z\right)  =-\left(  M-\delta\left\langle Mz,e_{1}\right\rangle
I\right)  z$ is the infinitesimal generator of the semigroup:%
\[
\varphi_{t}\left(  z\right)  =\frac{\exp\left(  tM\right)  z}{\delta
\left\langle z,\left(  \exp\left(  tM\right)  ^{H}-E\right)  e_{1}%
\right\rangle +1},
\]
according to Lemma \ref{infinity1}, $\left\{  \varphi_{t}\right\}  $ is a
semigroup of linear fractional self-maps of $B_{N},$ and%
\[
\varphi_{1}\left(  z\right)  =\varphi\left(  z\right)  .
\]

\end{proof}

\section{Non-elliptic cases}

Let $\varphi$ be a hyperbolic or a parabolic linear factional map. The
following lemma, which is a modified version of Theorem 4.1 of \cite{BCM1},
shows that $\varphi$ is conjugated to some special linear map on
$\mathbb{H}^{N}$. The concept of pseudo-inverse of a matrix is used in Theorem
\ref{semi4.1}. For more details about pseudo-inverse, we refer to \cite{Me}, p422.

\begin{theorem}
\label{semi4.1} Let $\varphi\in\operatorname{LFM}\left(  B_{N}\right)  $ be
non-elliptic with boundary dilation coefficient $\frac{1}{\lambda}$. Then
$\varphi$ is conjugated to a self-map $\psi$ of $\mathbb{H}^{N}$ which is
given by%
\[
\psi\left(  z,w\right)  =\left(  \lambda z+2i\left\langle w,a\right\rangle
+b,Mw+c\right)  \text{, \ }\left(  z,w\right)  \in\mathbb{H}^{N}%
\subset\mathbb{C}_{{}}\times\mathbb{C}^{N-1},
\]
where $c\in\mathbb{C}$, $b,d\in\mathbb{C}^{N-1}$, $M\in\mathbb{C}^{\left(
N-1\right)  \times\left(  N-1\right)  }$. Conversely, such a map is a self-map
of $\mathbb{H}^{N}$ if and only if

$\left(  P1\right)  \;Q:=\lambda I-M^{H}M$ is a Hermitian positive
semi-definite matrix;

$\left(  P2\right)  \;\operatorname{Im}\left(  b\right)  -\left\vert
c\right\vert ^{2}\geq\left\langle Q^{+}\left(  M^{\ast}c-a\right)  ,M^{\ast
}c-a\right\rangle $ where $Q^{+}$ is the pseudo-inverse of $Q;$

$\left(  P3\right)  \;QQ^{+}\left(  M^{\ast}c-a\right)  =M^{\ast}c-a$.
\end{theorem}

\begin{proof}
The only difference compared with Theorem 4.1 in \cite{BCM1} is $\left(
P3\right)  $. The corresponding condition there is $M^{\ast}c-a$ belongs to
the space spanned by the columns of $Q$. That is to say, there is a vector
$x\in\mathbb{C}^{N-1}$ such that%
\[
Qx=M^{\ast}c-a.
\]
According to the property of $Q^{+}$, the above equation has at least one
solution if and only if%
\[
QQ^{+}\left(  M^{\ast}c-a\right)  =M^{\ast}c-a.
\]

\end{proof}

\subsection{The parabolic cases}

The following theorem shows that a parabolic automorphism can always be
imbedded into a semigroup on $B_{N}$.

\begin{theorem}
\label{paraaut}Let $\varphi$ be a parabolic automorphism of $B_{N}$. Then

\begin{enumerate}
\item $\varphi$ is conjugated to $\psi\in\operatorname{Aut}\left(
\mathbb{H}^{N},\mathbb{H}^{N}\right)  $ which is defined by%
\[
\psi\left(  z,u,v\right)  =\left(  z+2i\left\langle u,a\right\rangle
+i\left\vert a\right\vert ^{2}+b,u+a,Dv\right)
\]
where $b$ is a real number, $a\in\mathbb{C}^{k}$ and $D\in\mathbb{C}^{\left(
N-k-1\right)  \times\left(  N-k-1\right)  }$ is diagonal, $\sigma\left(
D\right)  \subset\partial\mathbb{D}\backslash\left\{  1\right\}  $.

\item $\varphi$ can be embedded into a semigroup of $B_{N}$.
\end{enumerate}
\end{theorem}

\begin{proof}
\begin{enumerate}
\item The proof can be found in Proposition 4.3 in \cite{BCM1}

\item Since $0\notin\sigma\left(  D\right)  $, there is a diagonal matrix
$\Theta$ such that $\exp\left(  i\Theta\right)  =D$. Let%
\[
\psi_{t}\left(  z,u,v\right)  =\left(  z+2i\left\langle u,ta\right\rangle
+it^{2}\left\vert a\right\vert ^{2}+tb,u+ta,\exp\left(  it\Theta\right)
v\right)  ,
\]
then clearly, for every $t>0$, $\psi_{t}$ is an automorphism of $\mathbb{H}%
^{N}$ and $\left\{  \psi_{t}\right\}  $ is a semigroup of $\mathbb{H}^{N}$.
Therefore, $\psi$ can be embedded into some semigroup of $\mathbb{H}^{N}$.
Consequently, $\varphi$ can be embedded into a semigroup of $B_{N}$.
\end{enumerate}
\end{proof}

Now we turn to arbitrary parabolic linear fractional self-maps.

\begin{theorem}
\label{paraim}Let $\varphi\in\operatorname{LFM}\left(  B_{N}\right)  $ be
parabolic. Then

\begin{enumerate}
\item $\varphi$ is conjugated to $\psi:\mathbb{H}_{{}}^{N}\rightarrow
\mathbb{H}^{N}$ which is defined by%
\begin{equation}
\psi\left(  z,u,v,w\right)  =\left(  z+2i\left\langle u,a\right\rangle
+2i\left\langle w,c\right\rangle +b,u+a,Dv,Aw\right)  , \label{paranormalform}%
\end{equation}
where
\[
b\in\mathbb{C},a\in\mathbb{C}^{p},c\in\mathbb{C}^{q},D\in\mathbb{C}^{r\times
r},A\in\mathbb{C}^{\left(  N-p-q-r-1\right)  \times\left(  N-p-q-r-1\right)  }%
\]
with

\begin{enumerate}
\item $D$ is diagonal, $\sigma\left(  D\right)  \subset\partial\mathbb{D}%
\backslash\left\{  1\right\}  ;$

\item $Q=I-A^{H}A$ is a Hermitian positive semi-definite matrix;

\item $\operatorname{Im}\left(  b\right)  -\left\vert a\right\vert ^{2}%
\geq\left\langle Q^{+}c,c\right\rangle $;

\item $QQ^{+}c=c.$
\end{enumerate}

\item Suppose $0\notin\sigma\left(  A\right)  $. Let
\[
\exp\left(  M\right)  =A\text{,}%
\]
and%
\begin{align*}
c_{t}  &  =\left(  I-\exp\left(  M\right)  ^{H}\right)  ^{-1}\left(
I-\exp\left(  tM\right)  ^{H}\right)  c,\\
Q_{t}  &  =I-\exp\left(  tM\right)  ^{H}\exp\left(  tM\right)  .
\end{align*}
If

\begin{enumerate}
\item $M$ is dissipative;

\item for any $t\geq0$, $t\left(  \operatorname{Im}b-\left\vert a\right\vert
^{2}\right)  \geq\lambda^{-t}\left\langle Q_{t}^{+}c_{t},c_{t}\right\rangle ,$

\item $Q_{t}Q_{t}^{+}c_{t}=c_{t},$
\end{enumerate}
\end{enumerate}

then $\psi$ can be embedded into a semigroup of $\mathbb{H}^{N}$.
\end{theorem}

\begin{remark}
Sometimes there are no $u$ or no $v$ or no $w$ appeared in $\left(
\ref{paranormalform}\right)  $. And due to Theorem 4.4 of \cite{BCM1}, we may
assume that $a=0$ if $\varphi$ has at least one invariant slice.
\end{remark}

Before the proof of Theorem \ref{paraim}, we need the following easy lemma.

\begin{lemma}
For any $\alpha,\beta\in\mathbb{C}$, $a\in\mathbb{C}^{p},D\in\mathbb{C}%
^{q\times q},A\in\mathbb{C}^{r\times r}$, let%
\begin{align*}
\tau\left(  z,W\right)   &  =\left(  z+2i\left\langle u,a\right\rangle
+\beta,u+a,Dv,w\right)  ,\\
\rho\left(  z,W\right)   &  =\left(  z+2i\left\langle w,c\right\rangle
+\alpha,u,v,Aw\right)  .
\end{align*}
Then%
\begin{equation}
\tau\circ\rho=\rho\circ\tau. \label{commute}%
\end{equation}

\end{lemma}

\begin{proof}
[Proof of Theorem \ref{paraim}]%
\begin{enumerate}
\item See section 2 of \cite{B1}.

\item Let%
\begin{align*}
\tau_{\psi}\left(  z,W\right)   &  =\left(  z+2i\left\langle u,a\right\rangle
+\beta_{\psi},u+a,Dv,w\right)  ,\\
\rho_{\psi}\left(  z,W\right)   &  =\left(  z+2i\left\langle w,c\right\rangle
+\alpha_{\psi},u,v,Aw\right)  ,
\end{align*}
with $b=\alpha_{\psi}+\beta_{\psi}$ and $\operatorname{Im}\left(  \beta_{\psi
}\right)  =\left\vert a\right\vert ^{2},\alpha_{\psi}\in\mathbb{C}$.

Then $\operatorname{Im}\alpha_{\psi}=\operatorname{Im}b-\left\vert
a\right\vert ^{2}$ and%
\[
\tau_{\psi}\circ\rho_{\psi}=\rho_{\psi}\circ\tau_{\psi}=\psi\text{.}%
\]
Since $\tau_{\psi}$ is a parabolic automorphism, according to Theorem
\ref{paraaut}, $\tau_{\psi}$ can be embedded into the semigroup $\left\{
\tau_{\psi,t}\right\}  $ which is defined by%
\[
\tau_{\psi,t}\left(  z,u,v,w\right)  =\left(  z+2i\left\langle
u,ta\right\rangle +it^{2}\left\vert a\right\vert ^{2}+t\operatorname{Re}%
\beta_{\psi},u+ta,\exp\left(  t\Theta_{D}\right)  ,w\right)  .
\]
Let%
\[
\rho_{\psi,t}\left(  z,u,v,w\right)  =\left(  u+2i\left\langle w,c_{t}%
\right\rangle +t\alpha_{\psi},u,v,\exp\left(  tM\right)  w\right)  .
\]
Then $\rho_{\psi,t}$ is a self-map of $\mathbb{H}^{N}$ for every $t\geq0$
according to Theorem\ \ref{semi4.1}. When $t=0$, we have
\[
c_{0}=0,b_{0}=0,\exp\left(  0M\right)  =E.
\]
Thereby,%
\[
\rho_{\psi,0}\left(  z,u,v,w\right)  =\left(  z,u,v,w\right)  \text{.}%
\]
Direct computation shows that for any $s,t\geq0$,
\[
\rho_{\psi,s}\circ\rho_{\psi,t}=\rho_{\psi,t}\circ\rho_{\psi,s}=\rho
_{\psi,s+t}\text{.}%
\]
That $\rho_{\psi,t}$ converges uniformly on compact subset of $\ \mathbb{H}%
^{N}$ when $t\rightarrow0^{+}$ is clear. As a consequence, $\left\{
\rho_{\psi,t}\right\}  $ is a semigroup of $\mathbb{H}^{N}$. Let%
\[
\psi_{t}=\tau_{\psi,t}\circ\rho_{\psi,t},
\]
then easy computations show that%
\begin{align*}
\psi_{t+s}  &  =\tau_{\psi,t+s}\circ\rho_{\psi,t+s}\\
&  =\psi_{t}\circ\psi_{s}\\
&  =\psi_{s}\circ\psi_{t}.
\end{align*}
Therefore $\psi$ can be embedded into a semigroup of $\mathbb{H}^{N}.$
\end{enumerate}
\end{proof}

Before we prove Theorem \ref{paramain}, we need the following lemma.

\begin{lemma}
\label{inequality1}Let $a\in\mathbb{C}$ with $\operatorname{Re}a<0$. Then%
\[
\sup_{t>0}\frac{\left\vert 1-e^{at}\right\vert ^{2}}{t\left(  1-\left\vert
e^{at}\right\vert ^{2}\right)  }=\frac{\left\vert a\right\vert ^{2}%
}{-2\operatorname{Re}a}.
\]

\end{lemma}

\begin{proof}
According to basic integration formula and H\H{o}lder inequality, if $t\geq0$,
then
\[
\left\vert 1-e^{at}\right\vert ^{2}=\left\vert -a\int_{0}^{t}e^{au}%
du\right\vert ^{2}\leq t\left\vert a\right\vert ^{2}\int_{0}^{t}\left\vert
e^{au}\right\vert ^{2}du.
\]
And%
\begin{align*}
t\left(  1-\left\vert e^{at}\right\vert ^{2}\right)   &  =-2\operatorname{Re}%
a\cdot t\int_{0}^{t}e^{2u\operatorname{Re}a}du\\
&  =-2\operatorname{Re}a\cdot t\cdot\int_{0}^{t}\left\vert e^{au}\right\vert
^{2}du.
\end{align*}
As a conclusion,%
\[
\frac{\left\vert 1-e^{at}\right\vert ^{2}}{t\left(  1-\left\vert
e^{at}\right\vert ^{2}\right)  }\leq\frac{\left\vert a\right\vert ^{2}%
}{-2\operatorname{Re}a}.
\]
Since%
\[
\lim_{t\rightarrow0^{+}}\frac{\left\vert 1-e^{at}\right\vert ^{2}}{t\left(
1-\left\vert e^{at}\right\vert ^{2}\right)  }=\frac{\left\vert a\right\vert
^{2}}{-2\operatorname{Re}a},
\]
we conclude that%
\[
\sup_{t>0}\frac{\left\vert 1-e^{at}\right\vert ^{2}}{t\left(  1-\left\vert
e^{at}\right\vert ^{2}\right)  }=\frac{\left\vert a\right\vert ^{2}%
}{-2\operatorname{Re}a}.
\]

\end{proof}

\begin{proof}
[Proof of Theorem \ref{paramain}]Let%
\[
M=diag\left(  -u_{1}+iv_{1},\cdots,-u_{r}+iv_{r}\right)  .
\]
Then
\[
A=\exp\left(  M\right)  .
\]
Denote by%
\begin{align*}
c_{t}  &  =\left(  I-\exp\left(  M^{H}\right)  \right)  ^{-1}\left(
I-\exp\left(  tM^{H}\right)  \right)  c,\\
Q_{t}  &  =I-\left(  \exp\left(  tM\right)  \right)  ^{H}\exp\left(
tM\right)  ,\\
b_{t}  &  =tb.
\end{align*}
Since both $A$ and $M$ are normal matrices and $\left\Vert B\right\Vert
=\left\Vert \exp\left(  M\right)  \right\Vert $,
\[
\left\Vert \exp\left(  M\right)  \right\Vert \leq1.
\]
According to proposition \ref{dissi}, for any $t\geq0$,
\[
\left\Vert \exp\left(  tM\right)  \right\Vert \leq1.
\]
Therefore $Q_{t}$ is Hermitian positive semi-definite and
\[
Q_{t}^{+}=Q_{t}^{-1}=\left[
\begin{array}
[c]{cccc}%
\frac{1}{1-e^{-2tu_{1}}} &  &  & \\
& \frac{1}{1-e^{-2tu_{2}}} &  & \\
&  & \ddots & \\
&  &  & \frac{1}{1-e^{-2tu_{r}}}%
\end{array}
\right]  .
\]
Besides,%
\begin{align*}
c_{t}  &  =\left(  E-\exp\left(  M^{H}\right)  \right)  ^{-1}\left(
E-\exp\left(  tM^{H}\right)  \right)  c\\
&  =\left[
\begin{array}
[c]{cccc}%
\frac{1-e^{t\left(  -u_{1}-iv_{1}\right)  t}}{1-e^{\left(  -u_{1}%
-iv_{1}\right)  }} &  &  & \\
& \frac{1-e^{t\left(  -u_{2}-iv_{2}\right)  t}}{1-e^{\left(  -u_{2}%
-iv_{2}\right)  }} &  & \\
&  & \ddots & \\
&  &  & \frac{1-e^{t\left(  -u_{r}-iv_{r}\right)  t}}{1-e^{\left(
-u_{r}-iv_{r}\right)  }}%
\end{array}
\right]  c.
\end{align*}
As a result,%
\[
c_{t}^{H}Q_{t}^{+}c_{t}=c^{H}\Theta_{t}c,
\]
where%
\[%
\begin{array}
[c]{ll}%
\Theta_{t}=diag & \left(  \frac{\left\vert 1-e^{t\left(  -u_{1}+iv_{1}\right)
}\right\vert ^{2}}{\left\vert 1-\lambda_{1}\right\vert ^{2}\left(
1-e^{-2tu_{1}}\right)  },\cdots,\left.  \frac{\left\vert 1-e^{t\left(
-u_{r}+iv_{r}\right)  }\right\vert ^{2}}{\left\vert 1-\lambda_{r}\right\vert
^{2}\left(  1-e^{-2tu_{r}}\right)  }\right)  .\right.
\end{array}
\]
Denote by%
\[
b=\left(  \beta_{1},\beta_{2},\cdots,\beta_{p}\right)  ^{T}\text{,}%
\]
then%
\[
b_{t}^{H}Q_{t}^{+}b_{t}=\sum_{j=1}^{p}\frac{\left\vert 1-e^{t\left(
-u_{j}+iv_{j}\right)  }\right\vert ^{2}}{\left\vert 1-\lambda_{j}\right\vert
^{2}\left(  1-e^{-2tu_{j}}\right)  }\left\vert \beta_{j}\right\vert ^{2}.
\]
Let%
\[
g_{\lambda_{j}}\left(  t\right)  =\frac{1}{t}\frac{\left\vert 1-e^{t\left(
-u_{j}+iv_{j}\right)  }\right\vert ^{2}}{\left\vert 1-\lambda_{j}\right\vert
^{2}\left(  1-e^{-2tu_{j}}\right)  }.
\]
According to Lemma \ref{inequality1}, for $j=1,2,\cdots,r$,
\[
\sup_{t\geq0}g_{\lambda_{j}}\left(  t\right)  =\frac{1}{2u_{j}}\left(
u_{j}^{2}+v_{j}^{2}\right)  \frac{1}{\left\vert 1-\lambda_{j}\right\vert ^{2}%
}\text{.}%
\]
Furthermore,%
\begin{align*}
\sup_{t\geq0}\left\{  \frac{1}{t}c_{t}^{H}Q_{t}^{+}c_{t}\right\}   &
=\sum_{j=1}^{p}\frac{1}{2u_{j}}\left(  u_{j}^{2}+v_{j}^{2}\right)  \frac
{1}{\left\vert 1-\lambda_{j}\right\vert ^{2}}\left\vert \beta_{j}\right\vert
^{2}\\
&  =c^{H}\Theta c.
\end{align*}
Consequently%
\[
c_{t}^{H}Q_{t}^{+}c_{t}\leq t\operatorname{Re}c.
\]

Due to Theorem \ref{paraim}, $\psi$ can be embedded into a semigroup of
$\mathbb{H}^{N}.$
\end{proof}

\subsection{The hyperbolic case}

The following is a similar lemma as the first part of Theorem \ref{paraaut}.

\begin{lemma}
[\cite{BCM1}]\label{hhyper}Let $\varphi$ be a hyperbolic automorphism of
$B_{N}$, then $\varphi$ is conjugated to $\psi\in\operatorname{Aut}\left(
\mathbb{H}^{N},\mathbb{H}^{N}\right)  $ with%
\[
\psi\left(  z,W\right)  =\left(  \lambda z+b,\sqrt{\lambda}UW\right)
\]
where $\left(  z,W\right)  \in\mathbb{C}\mathbf{\times}\mathbb{C}^{N-1}$,
$b\in\mathbb{R},$ $\lambda>1$ and $U$ is a unitary matrix.
\end{lemma}

\begin{theorem}
\label{hypernormal}Let $\varphi\in\operatorname{LFM}\left(  B_{N}\right)  $ be
hyperbolic and has at least an invariant slice. Then

\begin{enumerate}
\item $\varphi$ is conjugated to $\psi:\mathbb{H}^{N}\rightarrow\mathbb{H}%
^{N}$ with%
\[
\psi\left(  z,u,v,w\right)  =\left(  \lambda z+b,\sqrt{\lambda}u,\sqrt
{\lambda}Dv,\sqrt{\lambda}Aw+c\right)  ,
\]
where $\lambda>1$ and

\begin{enumerate}
\item $D$ is diagonal, $\sigma\left(  D\right)  \subset\partial\mathbb{D}%
\backslash\left\{  1\right\}  ;$

\item Both $Q=I-A^{H}A$ and $P=I-AA^{H}$ are Hermitian positive semi-definite matrices;

\item $\operatorname{Im}\left(  b\right)  \geq\left\langle P^{+}%
c,c\right\rangle $;

\item $QQ^{+}A^{H}c=A^{H}c$.
\end{enumerate}

\item Suppose $A$ is non-singular and there exists a matrix $M$ such that%
\[
\exp\left(  M\right)  =B.
\]
Denoted by
\begin{align*}
\lambda_{t}  &  =\lambda^{t},\\
A_{t}  &  =\exp\left(  tM\right)  ,\\
a_{t}  &  =\left(  \lambda-\sqrt{\lambda}A^{H}\right)  ^{-1}\left(
\lambda_{t}-\sqrt{\lambda_{t}}A_{t}^{H}\right)  a,\\
b_{t}  &  =\frac{1-\lambda_{t}}{1-\lambda},\\
Q_{t}  &  =E-A_{t}^{H}A_{t}.
\end{align*}
If

\begin{enumerate}
\item $Q_{t}$ is Hermitian positive semi-definite;

\item for any $t\geq0$,
\[
\operatorname{Im}b_{t}\geq\frac{1}{\lambda_{t}}\left\langle Q_{t}^{+}%
a_{t},a_{t}\right\rangle ;
\]

\item for any $t\geq0$, $Q_{t}Q_{t}^{+}a_{t}=a_{t}$,
\end{enumerate}
\end{enumerate}

then $\psi$ can be embedded in to a semigroup of $\mathbb{H}^{N}$.
\end{theorem}

\begin{proof}
\begin{enumerate}
\item See \cite[Proposition 2.3]{hyperbolic}.

\item The proof of the above theorem is just the same with the second part of
Theorem \ref{paraim}, we omit it here.
\end{enumerate}
\end{proof}

The following corollary shows that a hyperbolic linear fractional map has
another form of normal form, for a proof, see \cite{BCM1}.

\begin{corollary}
Let $\varphi\in\operatorname{LFM}\left(  B_{N}\right)  $ be hyperbolic. Then
$\varphi$ is conjugated to $\psi\in\operatorname{LFM}\left(  \mathbb{H}%
^{N}\right)  $ with%
\[
\psi\left(  z,u,v,w\right)  =\left(  \lambda z+2i\left\langle w,a\right\rangle
+b,\sqrt{\lambda}u,\sqrt{\lambda}Dv,\sqrt{\lambda}Aw+c\right)  ,
\]
where $\lambda>1$ and

\begin{enumerate}
\item $D$ is diagonal, $\sigma\left(  D\right)  \subset\partial\mathbb{D}%
\backslash\left\{  1\right\}  ;$

\item $Q=I-A^{H}A$ is Hermitian positive semi-definite matrix;

\item $\operatorname{Im}\left(  b\right)  -\left\vert c\right\vert ^{2}%
\geq\left\langle Q^{+}\left(  A^{H}c-\frac{a}{\sqrt{\lambda}}\right)
,A^{H}c-\frac{a}{\sqrt{\lambda}}\right\rangle $;

\item $QQ^{+}\left(  A^{H}c-\frac{a}{\sqrt{\lambda}}\right)  =A^{H}c-\frac
{a}{\sqrt{\lambda}}$.
\end{enumerate}
\end{corollary}

\begin{lemma}
\label{equality2}Let $a\in\mathbb{C}$ with $\operatorname{Re}a<0$ and
$\lambda>1,$ $\lambda+2\operatorname{Re}a<0$. Then%
\[
\sup_{t>0}\frac{1}{1-e^{-\lambda t}}\cdot\frac{\left\vert 1-e^{at}\right\vert
^{2}}{1-e^{\lambda t}\left\vert e^{at}\right\vert ^{2}}=-\frac{\left\vert
a\right\vert ^{2}}{\lambda\left(  \lambda+2\operatorname{Re}a\right)  }.
\]

\end{lemma}

\begin{proof}
First of all,
\[
\lim_{t\rightarrow0^{+}}\frac{1}{1-e^{-\lambda t}}\cdot\frac{\left\vert
1-e^{at}\right\vert ^{2}}{1-e^{\lambda t}\left\vert e^{at}\right\vert ^{2}%
}=-\frac{\left\vert a\right\vert ^{2}}{\lambda\left(  \lambda
+2\operatorname{Re}a\right)  }.
\]
And%
\begin{align*}
\left\vert 1-e^{at}\right\vert ^{2}  &  =\left\vert a\right\vert
^{2}\left\vert \int_{0}^{t}e^{au}\cdot e^{\frac{\lambda}{2}u}\cdot
e^{-\frac{\lambda}{2}u}du\right\vert ^{2}\\
&  \leq\left\vert a\right\vert ^{2}\int_{0}^{t}\left\vert e^{\left(
\lambda+2a\right)  u}\right\vert ^{2}du\int_{0}^{t}e^{-\lambda u}du\\
&  =\left\vert a^{2}\right\vert \left(  \frac{e^{-\lambda t}-1}{-\lambda
}\right)  \frac{e^{\lambda t}\left\vert e^{at}\right\vert ^{2}-1}%
{\lambda+2\operatorname{Re}a}.
\end{align*}
Therefore%
\[
\frac{1}{1-e^{-\lambda t}}\cdot\frac{\left\vert 1-e^{at}\right\vert ^{2}%
}{1-e^{\lambda t}\left\vert e^{at}\right\vert ^{2}}\leq-\frac{\left\vert
a\right\vert ^{2}}{\lambda\left(  \lambda+2\operatorname{Re}a\right)  },
\]
and our lemma holds.
\end{proof}

\begin{proof}
[Proof of Theorem \ref{hypermain}]Let%
\[
M=diag\left(  -u_{1}+iv_{1},\cdots,-u_{r}+iv_{r}\right)  \text{,}%
\]
then%
\[
A=\exp\left(  M\right)  .
\]
Denote by%
\begin{align*}
\lambda_{t}  &  =\lambda^{t},\\
A_{t}  &  =\exp\left(  tM\right)  ,\\
a_{t}  &  =\left(  \lambda-\sqrt{\lambda}A^{H}\right)  ^{-1}\left(
\lambda_{t}-\sqrt{\lambda_{t}}A_{t}^{H}\right)  a,\\
b_{t}  &  =\frac{1-\lambda_{t}}{1-\lambda}b,\\
Q_{t}  &  =E-A_{t}^{H}A_{t}.
\end{align*}
Then%
\[
\left\langle Q_{t}^{+}a_{t},a_{t}\right\rangle =a^{H}diag\left(  \alpha
_{1}\left(  t\right)  ,\cdots,\alpha_{s}\left(  t\right)  \right)  a,
\]
where%
\[
\alpha_{j}\left(  t\right)  =\frac{\left\vert \lambda^{t}-\sqrt{\lambda^{t}%
}e^{t\left(  -u_{j}-iv_{j}\right)  }\right\vert ^{2}}{\left(  1-e^{-2tu_{j}%
}\right)  \left\vert \lambda-\sqrt{\lambda}e^{-u_{j}-tv_{j}}\right\vert ^{2}%
}=\frac{\lambda^{2t}\left\vert 1-e^{-t\left(  \frac{\ln\lambda}{2}%
+u_{j}\right)  +iv_{j}}\right\vert ^{2}}{\left(  1-e^{-2tu_{j}}\right)
\left\vert \lambda-\sqrt{\lambda}\lambda_{j}\right\vert ^{2}}.
\]
Notice that according to Lemma \ref{equality2}, for $\frac{\ln\lambda}%
{2}+u_{j}>0,v_{j}\geq0$ and $t\geq0$,%
\begin{align*}
\frac{\lambda^{t}\left\vert 1-e^{-t\left(  \frac{\ln\lambda}{2}+u_{j}\right)
+iv_{j}}\right\vert ^{2}}{\left(  1-e^{-2tu_{j}}\right)  \left(  \lambda
^{t}-1\right)  }  &  =\frac{\left\vert 1-e^{-t\left(  \frac{\ln\lambda}%
{2}+u_{j}\right)  +iv_{j}}\right\vert }{\left(  1-e^{t\ln\lambda}e^{-2t\left(
u_{j}+\frac{\ln\lambda}{2}\right)  }\right)  \left(  1-e^{-t\ln\lambda
}\right)  }\\
&  \leq\frac{1}{2u_{j}\ln\lambda}\left[  \left(  \frac{\ln\lambda}{2}%
+u_{j}\right)  ^{2}+v_{j}^{2}\right]  .
\end{align*}
Thus we get
\[
\frac{1}{\lambda^{t}}\frac{\left(  \lambda-1\right)  }{\left(  \lambda
^{t}-1\right)  }\left\langle Q_{t}^{+}a_{t},a_{t}\right\rangle \leq
a^{H}\Theta a\leq b.
\]
Our conclusion follows from Proposition \ref{hypernormal}.
\end{proof}

\section{The case of $N=2$ and case of automorphisms}

\begin{proof}
[Proof of Corollary \ref{autmain}]When $\varphi$ is an elliptic automorphism,
$\varphi$ is conjugated to a unitary transformation of $B_{N}$ and therefore
$\varphi$ can always be embedded into a semigroup of $B_{N}$.

Theorem \ref{paraaut} shows that a parabolic automorphism is always embeddable.

If $\varphi$ is a hyperbolic automorphism, then by Lemma \ref{hhyper},
$\varphi$ is conjugated to%
\[
\psi_{1}\left(  z^{\prime},z^{\prime\prime}\right)  =\frac{1}{\alpha}\left(
z^{\prime}+ic,\sqrt{\alpha}Uz^{\prime\prime}\right)  ,
\]
where $U$ is a unitary matrix, thus by Theorem \ref{hypermain}, $\varphi$ can
be embedded into a semigroup of $B_{N}$.
\end{proof}

\begin{proof}
[Proof of Theorem \ref{main3}]According to Theorem \ref{paranormalform}, it is
easy to see that $\left(  1\right)  $ holds.

Since $\psi_{2}$ and $\psi_{3}$ are all automorphisms, $\psi_{2}$ and
$\psi_{3}$ can always be embedded into some semigroups. $\psi_{1}$ is
embeddable follows from theorem \ref{paramain}.
\end{proof}

\begin{proof}
[Proof of Theorem \ref{main4}]%
\begin{enumerate}
\item According to Theorem \ref{hypernormal}, $\varphi$ is conjugated to%
\[
\tilde{\psi}\left(  u_{1},u_{2}\right)  =\left(  \lambda u_{1}+2i\left\langle
u_{2},\tilde{b}\right\rangle +\tilde{c},\mu u_{2}+\tilde{d}\right)  .
\]

Firstly, we assume that $\mu\neq1$. Let%
\[
\phi\left(  u_{1},u_{2}\right)  =\frac{1}{\beta}\left(  u_{1}+\frac{2}%
{\sqrt{\beta}}\bar{d}u_{2}+e,\sqrt{\beta}u_{2}+d\right)
\]
with $\operatorname{Re}e=\left\vert d\right\vert ^{2}$. Then $\phi$ is an
automorphism of $\mathbb{H}^{2}$ and%
\[
\phi^{-1}\left(  u_{1},u_{2}\right)  =\left(  -\frac{1}{\beta}\left(  \beta
e-u_{1}\beta^{2}-2\left\vert d\right\vert ^{2}+2\bar{d}u_{2}\beta\right)
,-\frac{1}{\sqrt{\beta}}\left(  d-u_{2}\beta\right)  \right)  .
\]
Now%
\[
\phi^{-1}\circ\tilde{\psi}\circ\phi\left(  u_{1},u_{2}\right)  =\left(
\ast,\frac{1}{\sqrt{\beta}}\left(  c\beta-d+d\mu+\sqrt{\beta}\mu u_{2}\right)
\right)  .
\]
Since $\mu\neq1$, let $d=\frac{c\beta}{1-\mu}$, then there exists $b$ and $c$
and $\alpha$ such that $\tilde{\psi}$ is conjugated to%
\[
\psi_{1}\left(  u_{1},u_{2}\right)  =\left(  \lambda u_{1}+2i\left\langle
u_{2},b\right\rangle +c,\sqrt{\lambda}\alpha u_{2}\right)  .
\]

Next, if $\mu=1$, then according to Theorem \ref{hypernormal}, $\psi$ is
conjugated to
\[
\psi_{2}\left(  u_{1},u_{2}\right)  =\left(  \lambda u_{1}+a,u_{2}+b\right)
.
\]
As an conclusion, $\left(  1\right)  $ holds.

\item $\left(  2\right)  $ holds following Theorem \ref{hypermain}.

\item Since
\[
\frac{d}{dt}\left(  \frac{\lambda^{t}-1}{t}\right)  =\frac{1}{t^{2}}\left(
t\lambda^{t}\ln\lambda-\lambda^{t}+1\right)
\]
and%
\[
\frac{d}{dt}\left(  t\lambda^{t}\ln\lambda-\lambda^{t}+1\right)  =t\lambda
^{t}\ln^{2}\lambda,
\]
thus for any $t>0$, we obtain%
\[
\frac{\lambda^{t}-1}{t}\geq\lim_{t\rightarrow0}\frac{\lambda^{t}-1}{t}%
=\ln\lambda.
\]
As a consequence,
\[
\frac{t^{2}}{\left(  \lambda^{t}-1\right)  ^{2}}\leq\frac{1}{\ln^{2}\lambda}.
\]
Let
\[
\tau_{t}\left(  u_{1},u_{2}\right)  =\left(  \lambda^{t}u_{1}+\frac
{\lambda^{t}-1}{\lambda-1}a,u_{2}+tb\right)  \text{,}%
\]
For any $t\geq0$,
\begin{align*}
\operatorname{Im}\left(  \frac{\lambda^{t}-1}{\lambda-1}a\right)   &
\geq\frac{\lambda^{t}-1}{\lambda-1}\cdot\frac{\left(  \lambda-1\right)  }%
{\ln^{2}\lambda}\left\vert b\right\vert ^{2}\\
&  \geq\left\vert tb\right\vert ^{2}.
\end{align*}
According to Theorem \ref{hypernormal}, $\tau_{t}$ is a self-map of
$\mathbb{H}^{2}$ and thus $\left\{  \tau_{t}\right\}  $ is a semigroup on
$\mathbb{H}^{2}$. Hence $\psi_{2}$ can be embedded into a semigroup of
$\mathbb{H}^{2}$.
\end{enumerate}
\end{proof}

\end{document}